\documentclass[graybox]{svmult}
\usepackage{mathptmx}       
\usepackage{helvet}         
\usepackage{courier}        
\usepackage{type1cm}        
\usepackage{makeidx}         
\usepackage{graphicx}        
\usepackage{multicol}        
\usepackage[bottom]{footmisc}

\usepackage{amscd}
\usepackage{amsfonts}
\usepackage{amsmath}
\usepackage{amssymb}
\usepackage{amstext}
\usepackage{dsfont}
\usepackage{fancyhdr}
\usepackage{graphicx}
\usepackage[latin1]{inputenc}
\usepackage{mathrsfs}
\usepackage{mathtools}
\usepackage{psfrag}
\usepackage{setspace}
\usepackage[normalem]{ulem}
\usepackage{verbatim}
\usepackage{url}

\usepackage[utf8]{inputenx}
\usepackage[T1]{fontenc}

\usepackage[hidelinks,colorlinks=true,urlcolor=blue,linkcolor=black,citecolor=black]{hyperref}
\newcommand*{\doi}[1]{\href{http://dx.doi.org/\detokenize{#1}}{doi}}

\allowdisplaybreaks

\renewcommand{\E}{\mathbb{E}}
\newcommand{\N}{\mathbb{N}}
\newcommand{\Pb}{\mathbb{P}}

\newcommand{\Z}{\mathbb{Z}}

\newcommand{\sleep}{\mathfrak{s}}
\newcommand{\toppling}{\mathfrak{t}}
\newcommand{\jump}{\mathfrak{m}}

\newcommand{\T}{T}

\newcommand{\CC}{\mathscr{C}}

\newcommand{\vep}{\varepsilon}

\newcommand{\oo}{0}

\newcommand{\dd}{{\mathrm d}}

\renewcommand{\leq}{\leqslant}
\renewcommand{\le}{\leqslant}
\renewcommand{\geq}{\geqslant}
\renewcommand{\ge}{\geqslant}
\renewcommand{\epsilon}{\varepsilon}

\makeindex             

\begin{document}

\title*{Avalanches in Critical Activated Random Walks}
\author{Manuel Cabezas and Leonardo T. Rolla}
\institute{Manuel Cabezas \at PUC-Chile, \email{mncabeza@mat.puc.cl}
\and Leonardo T. Rolla \at IMAS-Conicet and NYU-Shanghai \email{leorolla@dm.uba.ar}}
%
%
\maketitle

\abstract*{We consider Activated Random Walks on $\Z$ with totally asymmetric jumps and critical particle density, with different time scales for the progressive release of particles and the dissipation dynamics.
We show that the cumulative flow of particles through the origin rescales to a pure-jump self-similar process which we describe explicitly.}

\abstract{We consider Activated Random Walks on $\Z$ with totally asymmetric jumps and critical particle density, with different time scales for the progressive release of particles and the dissipation dynamics.
We show that the cumulative flow of particles through the origin rescales to a pure-jump self-similar process which we describe explicitly.}

\bigskip

\\


\begin{figure}[b]
\centering
\includegraphics[width=.8\textwidth]{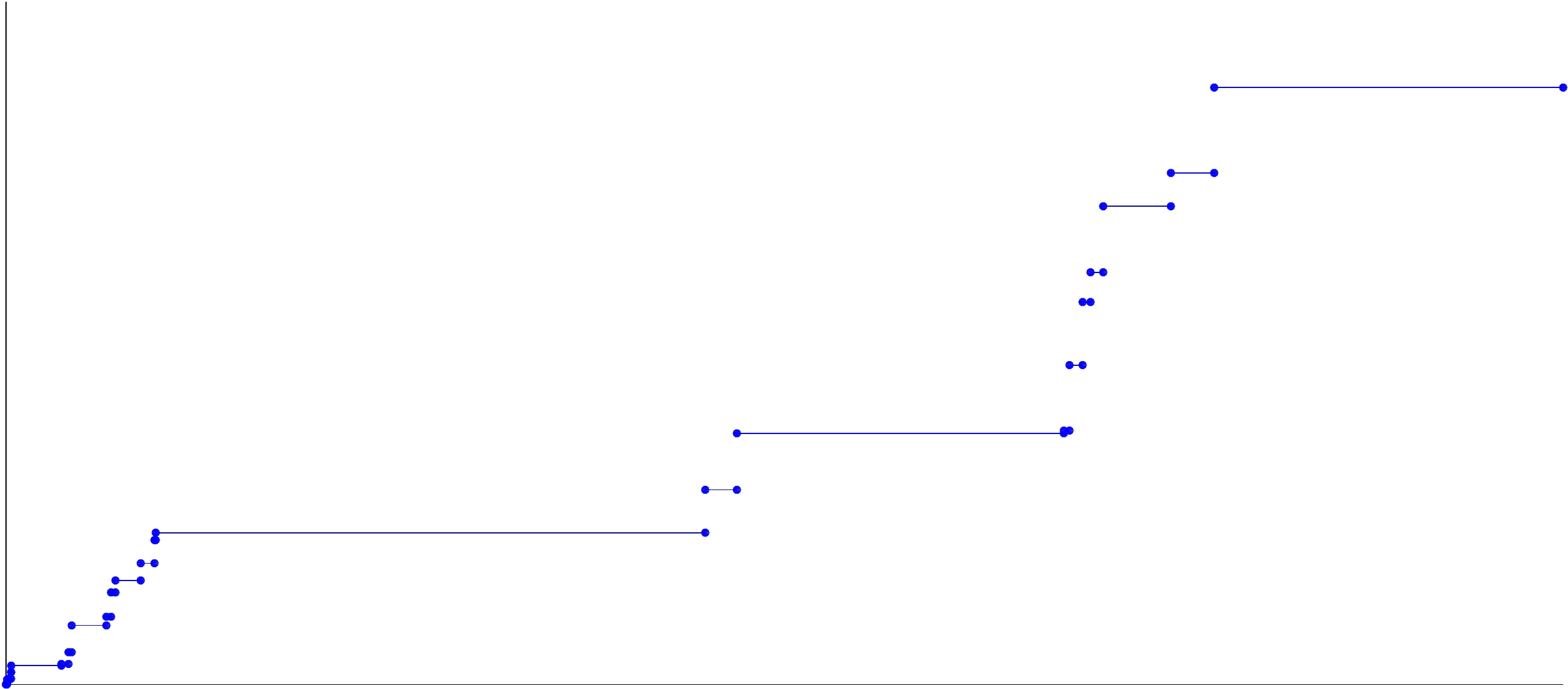}
\caption{\small Simulation of the avalanche process with $n=10^5$, $\zeta=.808$ and $\rho=.1$}
\end{figure}

\section{Introduction}
\label{sec:intro}

The totally asymmetric Activated Random Walk (ARW) dynamics on $\Z$ is a continuous-time conservative system made of active and passive particles, where each active particle jumps from $x$ to $x+1$ at rate $1$, and spontaneously decays to a passive state at rate $0 < \lambda < \infty$. Active particles reactivate passive particles instantly when they occupy the same site, in particular active particles at the same site prevent each other from decaying. This model has received increasing attention~\cite{AmirGurel-Gurevich10, AsselahSchapiraRolla19,BasuGangulyHoffman18, BasuGangulyHoffmanRichey19, CabezasRollaSidoravicius14, DickmanRollaSidoravicius10, RollaSidoravicius12, RollaTournier18, Shellef10, SidoraviciusTeixeira17, StaufferTaggi18, Taggi16, Taggi19}, see~\cite{Rolla19} for a self-contained introduction.

This model displays a phase transition in terms of the density of particles: Let $\eta(x)$ denote the initial number of particles at $x\in\Z$ and assume that the initial configuration is i.i.d.\ with mean $\zeta = \E [ \eta(\oo) ]$.
If $\zeta>\frac{\lambda}{1+\lambda}$ then the system can sustain a non-vanishing density of active particles, whereas, if $\zeta\leq\frac{\lambda}{1+\lambda}$, the density of active particles decays to $0$, see~\cite{Rolla19}. 
 
In this paper we study the flow process defined as follows.
Let $\eta_0 = \eta(0) \delta_0$, that is, the configuration where site $0$ has $\eta(0)$ particles and other sites are vacant.
We then run the above dynamics starting from $\eta_0$ until we get a configuration $\eta_0'$ without active particles, which we call \emph{stable}.
Finally, we let $C_0$ denote the number of particles which jump from $0$ to $+1$ during this evolution.
In the next step, we take $\eta_1 = \eta_0' + \eta(-1)\delta_{-1}$, that is, we add $\eta(-1)$ particles to the site $-1$.
Again we run the above dynamics starting from $\eta_1$ until we get a stable configuration $\eta_1'$, and let $C_1$ denote the total number of particles which jump from $0$ to $+1$.
In the same fashion, we define $\eta_{n}=\eta_{n-1}' + \eta(-n)\delta_{-n}$, stabilize $\eta_n$ obtaining $\eta_n'$, and define $C_n$ as the number of particles of $\eta_n'$ found to the right of site $0$.
Finally, the \emph{flow process} is defined as $(C_n)_{n=0,1,2,\dots}$.

Note that the sequence $(\eta_n')_{n}$ is a non-homogeneous Markov process with respect to its natural filtration. The transition probabilities are determined by the common distribution of $\eta(x)$ and by the dynamics run between each pair of steps, which is parametrized by $\lambda$.
Also, $C_n$ can be read from $\eta_n'$, but $(C_n)_n$ is not Markovian.



If the system is subcritical (i.e., $\zeta < \frac{\lambda}{1+\lambda}$) then $C_n$ is eventually constant.
It it is supercritical, $\frac{C_n}{n}$ tends to a positive number.
In the critical case, none of the above happens: $C_n$ diverges but $\frac{C_n}{n}$ vanishes, and one expects the system to have a non-trivial scaling limit.

Let $\sigma_s^2=\zeta-\zeta^2$ denote the variance of a Bernoulli variable with parameter $\zeta$, and let $\sigma_p^2 = \E[\eta(\oo)^2]-\zeta^2 \geqslant \sigma_s^2$.
Consider the critical case $\zeta = \frac{\lambda}{1+\lambda}$, assume that $\sigma_p^2 < \infty$ and define $\rho:=\frac{\sigma_s}{\sigma_p}\in(0,1]$.

\begin{figure}[t]
\centering
\includegraphics[width=.47\textwidth]{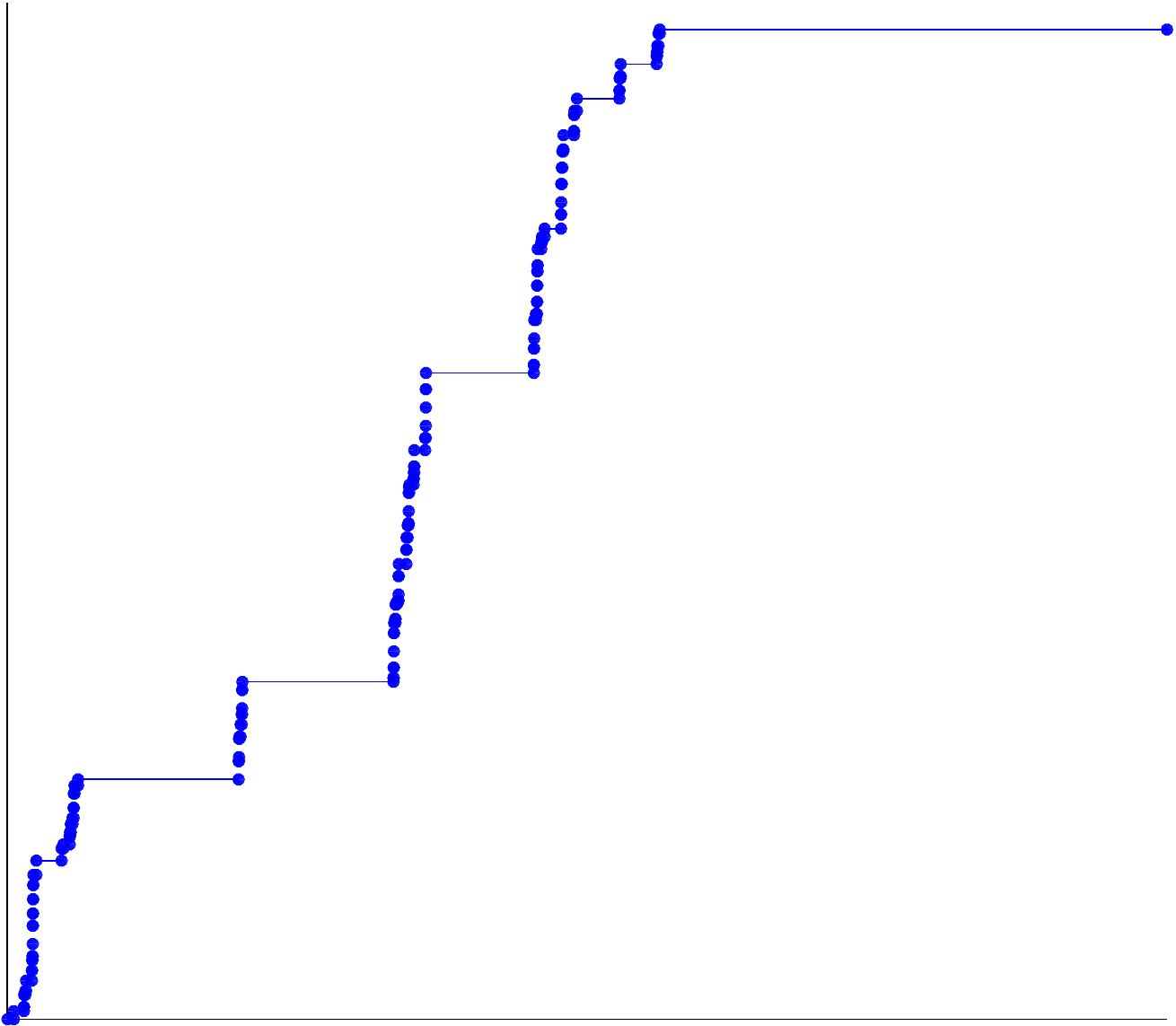}
\hfill{}
\includegraphics[width=.47\textwidth]{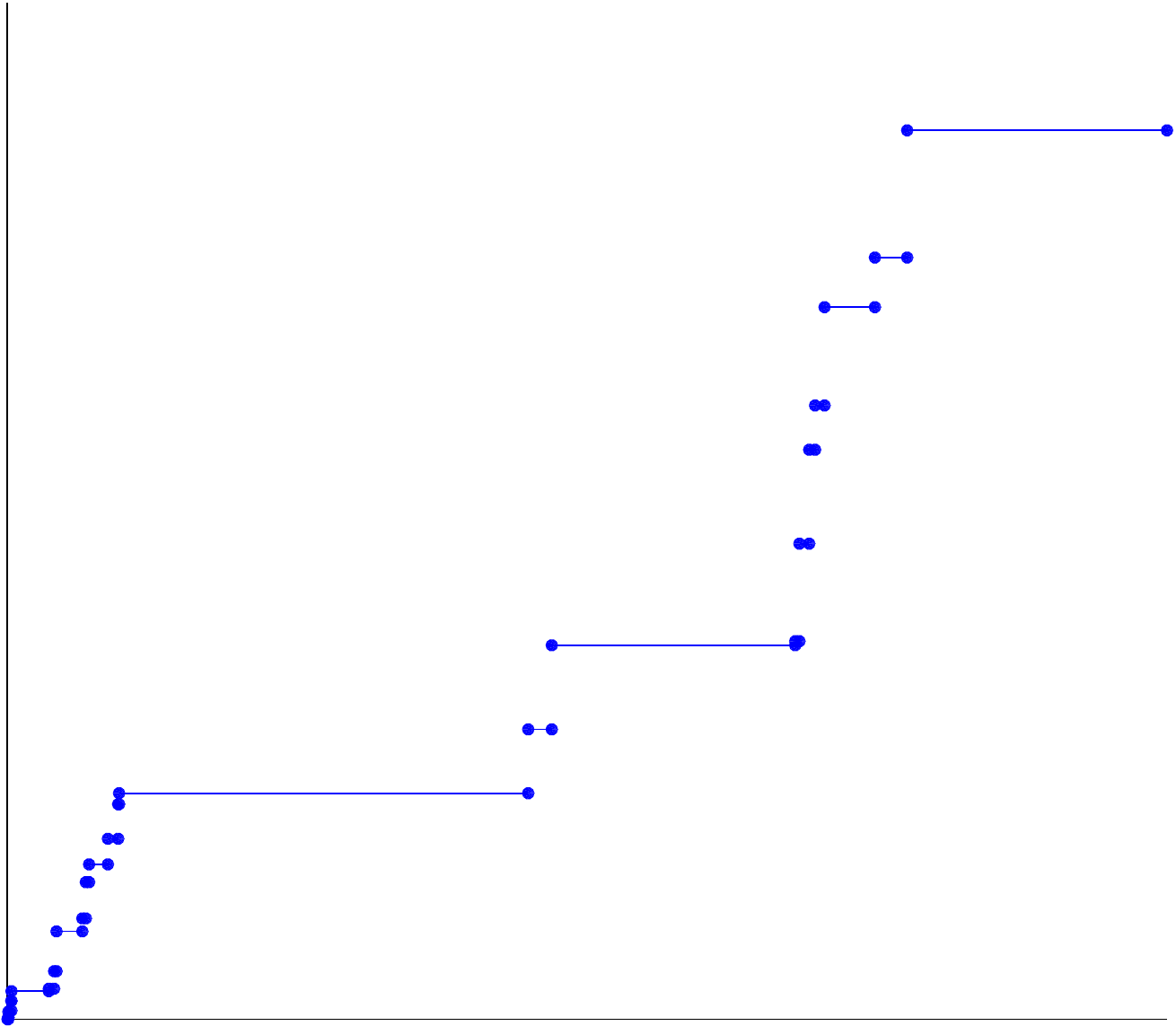}
\hfill{}
\caption{comparison of an avalanche process with $\rho=0$ and $\rho=0.1$}
\label{fig:newresult}
\end{figure}

\begin{theorem}
\label{thm:convergence}There exists a stochastic process $(\CC^\rho_x)_{x\ge 0}$ such that
\[
\frac{1}{\sigma_p}
\left(\varepsilon \, C_{\lfloor \varepsilon^{-2}x \rfloor}\right)_{x \ge 0}
\stackrel{\dd}{\to}
\left( \CC^\rho_{x}\right)_{x \ge 0}
\]
in the Skorohod $J_1$ metric as $\vep \to 0$.
\end{theorem}

\begin{remark}
The study of critical flows for the ARW started in~\cite{CabezasRollaSidoravicius14}, where the extreme case $\lambda = \infty$ was studied.
In this case, $\rho = 0$ and the limiting process $\left( \CC^0_{x}\right)_{x \ge 0}$ is the running maximum of a Brownian motion.%
\footnote
{In this case, the same convergence as stated in Theorem~\ref{thm:convergence} follows from the more complicated, continuous-time analysis done in~\cite{CabezasRollaSidoravicius14}, or alternatively from the arguments presented here.}
See Figure~\ref{fig:newresult} for a comparison.
\end{remark}

\begin{figure}[t]
\centering
\includegraphics[width=.31\textwidth]{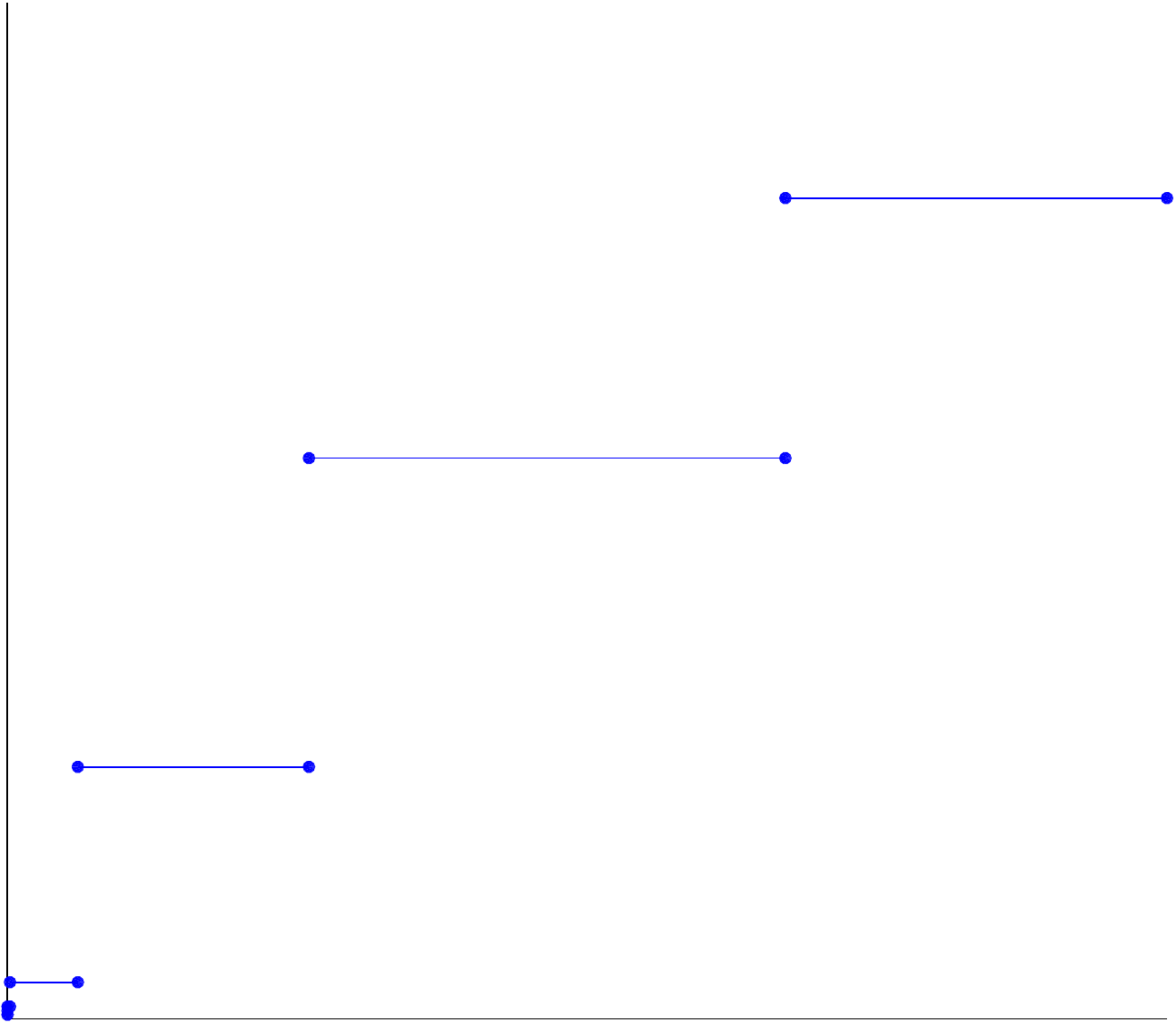}
\hfill
\includegraphics[width=.31\textwidth]{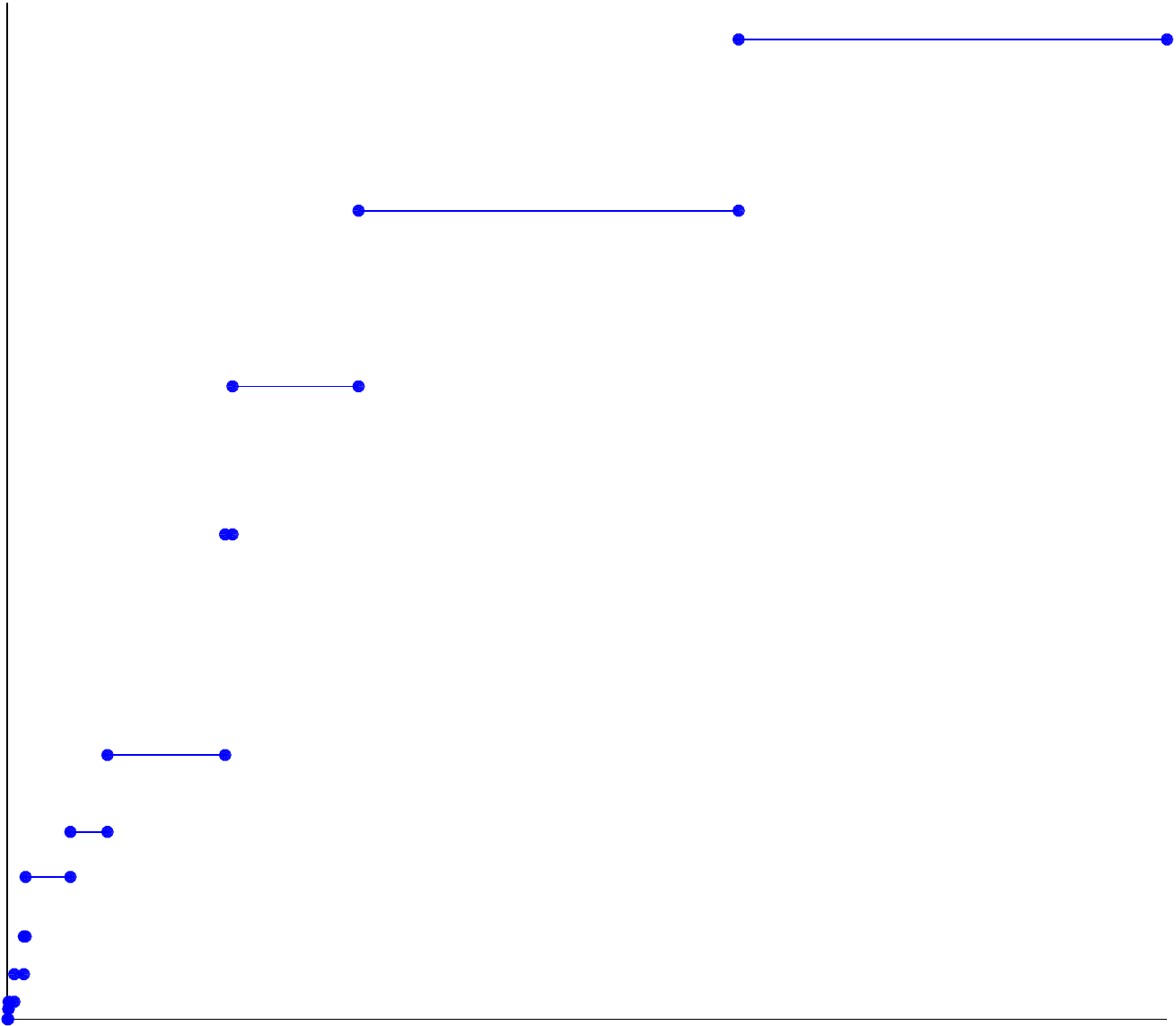}
\hfill
\includegraphics[width=.31\textwidth]{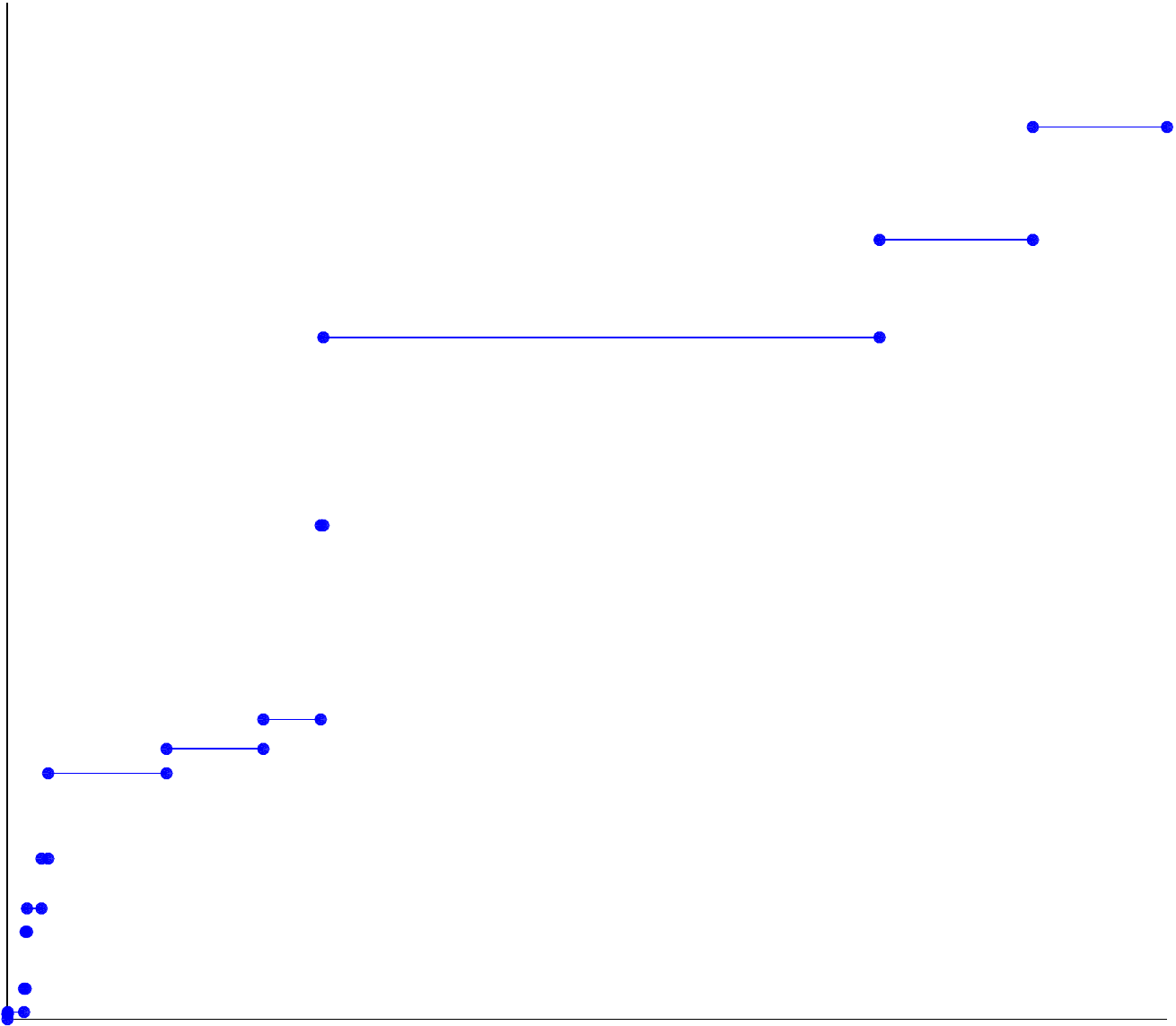}
\\[1em]
\includegraphics[width=.31\textwidth]{figures/run-0100-0808-20-11-100000-crop.pdf}
\hfill
\includegraphics[width=.31\textwidth]{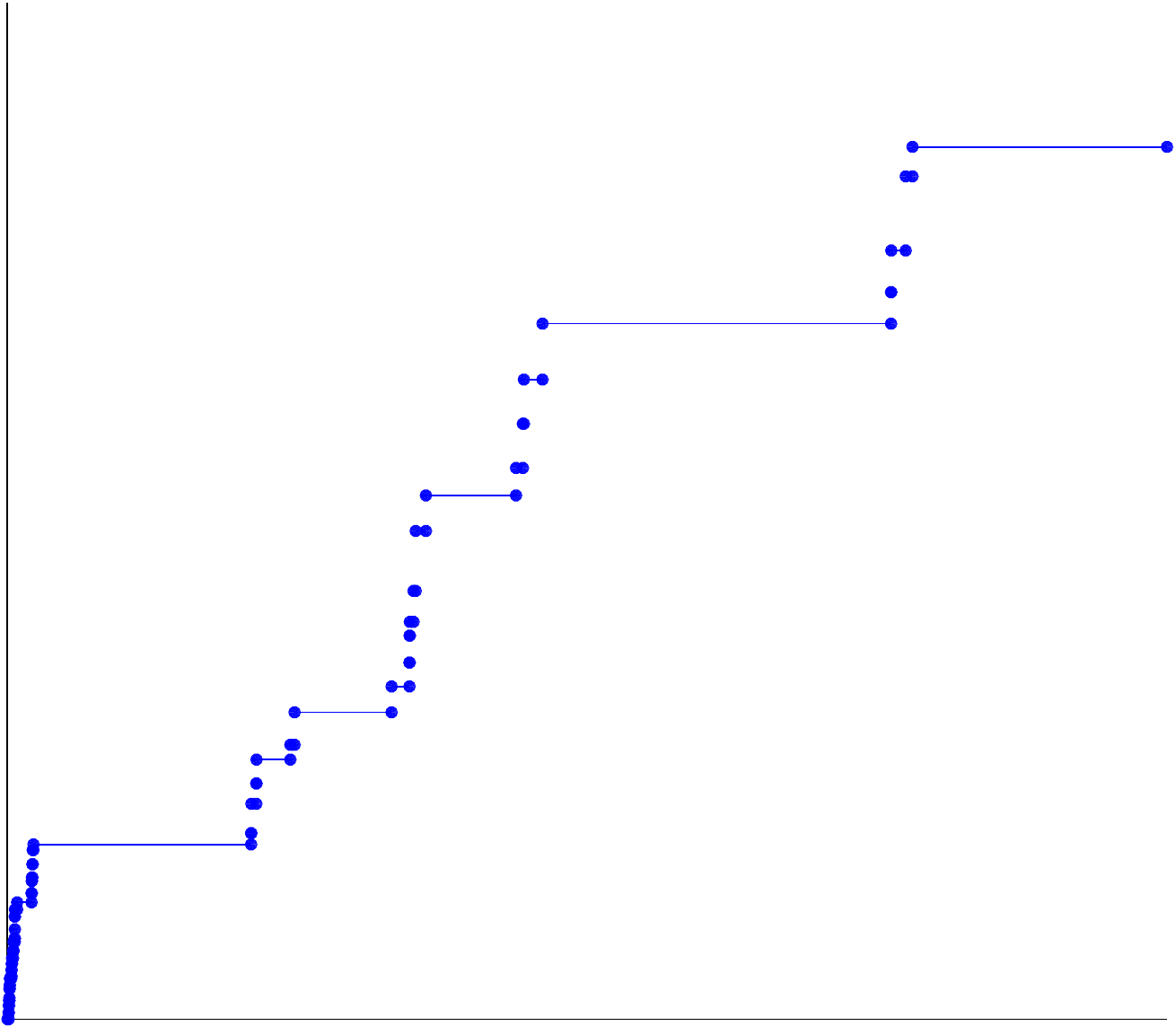}
\hfill
\includegraphics[width=.31\textwidth]{figures/run-0000-1000-10-08-25000-crop.pdf}
\caption{Simulations of $(C_n)_n$ with $\rho=1.000,\, 0.500,\, 0.300,\, 0.100,\, 0.051,\, 0.000$. The steps are $100000$ for the first four graphs, then 50000 and 25000.
For large $\rho$, the process seems more jumpy, and as $\rho$ gets smaller the size of jumps become smaller and jumps tend do cluster together, until finally at the extreme case $\rho=0$ studied in~\cite{CabezasRollaSidoravicius14} the process becomes continuous.}
\label{fig:manyrho}
\end{figure}

Below we will give a description of the stochastic process $\left( \CC^\rho_{x}\right)_{x \ge 0}$ in terms of a family of coalescing reflected correlated Brownian motions.
We point out that the superscript $\rho$ in the limiting process cannot be reduced to a multiplicative factor, and these processes are indeed very different as $\rho$ varies, see Figure~\ref{fig:manyrho}.

Before that, we state a qualitative property of the scaling limit which explains why we call it an avalanche process.
Being a scaling limit, it is scale invariant.

\begin{theorem}
\label{thm:purejump}
The process $\left( \CC^\rho_{x}\right)_{x \ge 0}$ consists of pure jumps.
Its jump times accumulate at 0 (as $ \CC^\rho $ is scale-invariant) but are otherwise discrete.
\end{theorem}

The limiting process $\CC^\rho$ can be informally described as follows.
Suppose $0 < x_0 < \cdots < x_k$, we want to sample $\CC^\rho_x$ for $x=x_0,\dots,x_k$.
Let $(\mathscr{P}_{x})_{x \leq 0}$ be a backward Brownian motion, with diffusivity constant equal to $1$, started at the origin.
Then, for each $i=0,\dots,k$, we run a Brownian motion $(\mathscr{R}^{i}_t)_{t\geq -x_i}$ started at $\mathscr{R}^{i}_{-x_i}=\mathscr{P}_{-x_i}$, with diffusivity $\rho$, reflected from below by the graph of $\mathscr{P}$.
We let the $k+1$ reflected Brownian motions diffuse independently (except for the reflection) until they coalesce.
Finally, $\CC_{x_i}$ is given by ${\mathscr{R}}^{i}_{0}$ for $i=0,\dots,k$.

Finally, we would like to say a word on the metric of Theorem~\ref{thm:convergence}. The convergence being in the Skorohod $J_1$-metric implies that each jump in the limiting process must be matched by a corresponding jump in the discrete process.
Hence, the discontinuities in $\CC^\rho$ (guaranteed by Theorem~\ref{thm:purejump}) correspond to abrupt increments, avalanches in the discrete process $(C_i)_{i\in\N}$.

The entire literature on ARW would probably not exist were it not for our dear friend Vladas Sidoravicius.
He was always in contact with prominent physicists, including Ronald Dickman, bringing interesting cutting-edge problems to the Probability community and enthusiastically promoting them.
In particular, he has been advertising the Activated Random Walk model since the early 2000s. In 2007, he proposed this problem to the second author as part of his PhD studies, which finally resulted in~\cite{RollaSidoravicius12}.
During our joyful meetings in the early 2010s, we worked on predecessors of the current paper~\cite{CabezasRollaSidoravicius14,CabezasRollaSidoravicius18}.

This paper is organized as follows.
In \S\ref{section:limit} we give the formal definition of the limiting process $\CC^\rho$.
In \S\ref{sec:construction} we define the ARW and state the Abelian property of its site-wise construction.
In \S\ref{section:seq} we introduce the sequential stabilization that will bet used throughout the article.
In \S\ref{sec:convergence} we prove convergence of finite dimensional distributions to later get full convergence (in the $J_1$ Skorohod metric) in \S\ref{sec:tightness}.
Finally, in \S\ref{sec:purejump} we show that $\CC^\rho$ is a pure-jump process when $\rho\neq0$.

\section{Formal definition of the limiting process}
\label{section:limit}


Next, we will describe formally the finite-dimensional distributions of the limiting process.
For convenience, and building the connection with the upcoming proofs, we describe the differences $\mathscr{R}^i - \mathscr{P}$, instead, given by the $B^{+,i}$ below.


Let $k\in\N_0$ and $-x_k\leq -x_{k-1} \leq \dots< -x_0 \leq0$.
Let $(\mathscr{S}^i_x)_{x\geq-x_i}, i=0,\dots,k.$ be independent Brownian motions started at $0$, $(\mathscr{S}^i_{-x_i}=0)$ with diffusion coefficient $\sigma_s$. Let also $(\mathscr{P}_x)_{x\leq0}$ be a backwards Brownian motion started at $0$ ($\mathscr{P}_{0}=0$) with diffusion coefficient $\sigma_p$ and independent of $\mathscr{S}^i,i=0,\dots,k$. For $i=0,\dots,k$ and $x\geq-x_i$ define
\[\tilde{B}^i_x:=\mathscr{P}_x-\mathscr{P}_{-x_i}-\mathscr{S}^i_x,\]
so that $(\tilde{B}^i_x)_{x\geq-x_i}$ are Brownian motions started at $0$, that is $(\tilde{B}^i_{-x_i}=0)$, with diffusion coefficient $r=\sqrt{\sigma_s^2+\sigma_p^2}$.
For $i=0,\dots,k$ and $x\geq-x_i$ we define
\begin{equation}
\label{eq:reflect}
\tilde{B}^{+,i}_x:=\tilde{B}^i_x-\inf\{\tilde{B}^i_s:s\in[-x_i,x]\}
\end{equation}
so that $(\tilde{B}^{+,i}_x)_{x\geq-x_i}, i=0,\dots,k$ are reflected Brownian motions started at $0$ with diffusion coefficient $r$.
Let $B^{+,k}:=\tilde{B}^{+,k}$ and, for $i=k-1,k-2,\dots,0$, let
\begin{equation}
\label{eq:defcoalesrefectbrown}
\tau_i:=\inf\{x\geq-x_i: \tilde{B}^{+,i}_x={B}^{+,i+1}_x\}
\quad
\text{ and }
\quad
{B}^{+,i}_x:=\begin{cases}
\tilde{B}^{+,i}_x &: x \in[-x_i,\tau_i),\\
{B}^{+,i+1}_x &: x \geq \tau_i.
\end{cases}
\end{equation}

Then, for each $k\in\N$ and each sequence $0\leq x_1\leq\cdots \leq x_k$ we have that
\[
(\CC^\rho_{x_0},\dots,\CC^\rho_{x_k})\stackrel{d}{=}(\tfrac{1}{\sigma_p}B^{+,0}_0,\dots,\tfrac{1}{\sigma_p} B^{+,k}_0)
,
\]
where $\stackrel{d}{=}$ denotes equality in distribution.

The above description is self-consistent, in the sense that removing points from the $\{x_0,\dots,x_k\}$ does not affect the distribution of the remaining points. Since the resulting process is non-decreasing, we can already deduce the existence of a càdlàg process $\left( \CC^\rho_{x}\right)_{x \ge 0}$ with these finite-dimensional distributions.

\begin{figure}[t]
\centering
\includegraphics[width=\textwidth]{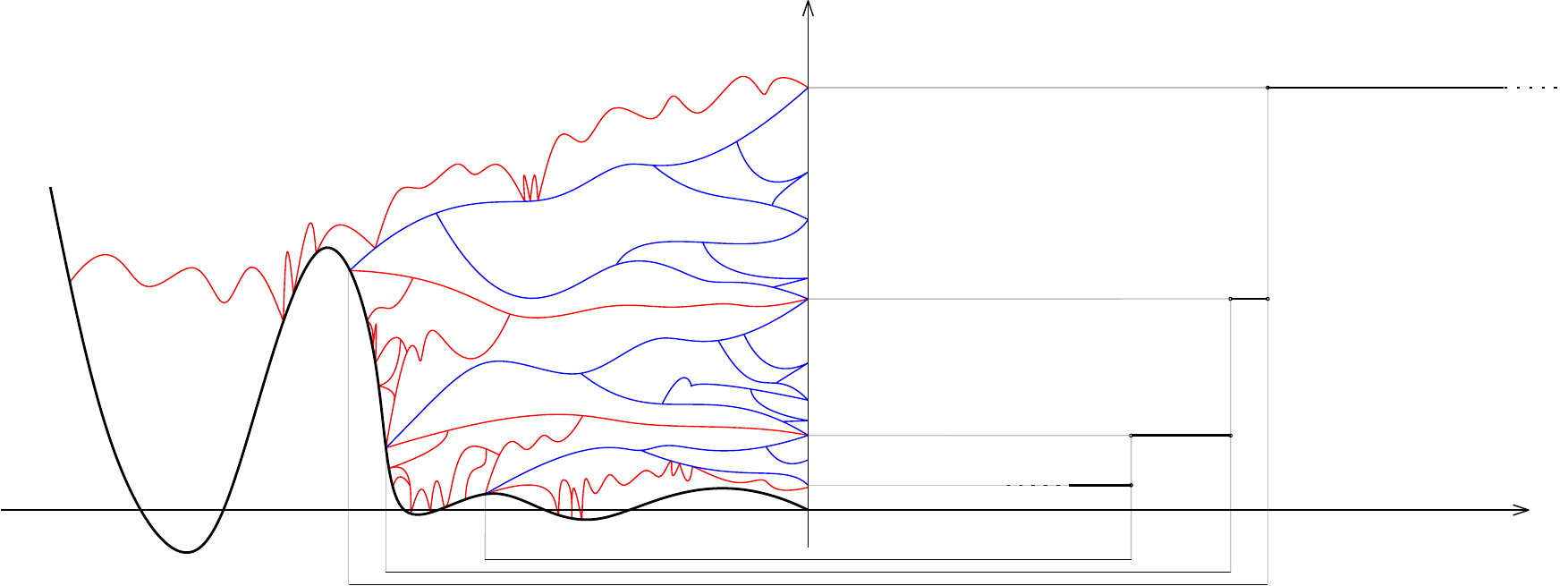}
\caption{Black path, blue paths and dual red paths on the left, $\CC^{\rho}$ on the right.}
\label{fig:artistic}
\end{figure}

In Figure~\ref{fig:artistic}, we see the graph of $\CC^\rho$ on the right-hand side ($x\geq0$) and the elements of this construction on the left side.
There is a black curve starting at $(0,0)$ moving backwards (from right to left). This curve is the Brownian motion $\mathscr{P}$ with diffusivity $\sigma_p$. 
From each point on the black we can start a red path moving forward (from left to right). This red path is a Brownian with diffusivity $\sigma_s$, reflected against the black curve. Different red paths diffuse independently until they meet, and coalesce after that.
To find the value of $\CC^{\rho}$ at a certain $x_0>0$, we follow the red path which starts at $(-x_0,\mathscr{P}_{-x_0})$ until it hits the vertical axis $\{x=0\}$.
The value $\CC^{\rho}_{x_0}$ is the height coordinate of the red path at this terminal point.

An alternative description is the following.
By considering the dual of the red paths one obtains a family of blue paths.
These paths start from each point on the positive vertical axis and move from right to left.
Blue paths coalesce when they meet, and they terminate upon hitting the black curve.
Each terminal point determines an interval given by the initial height of the blue paths which terminate there.
This interval corresponds to the jumps of $\CC^\rho$, see Figure~\ref{fig:artistic}.
On the left side we see the intervals determined by the blue paths. The extremes of the interval (down and up) correspond to the different (left and right) limits at a discontinuity of $\CC^\rho$. The $x$ coordinate of the jumps is (the reflection of) the $x$ coordinate of the terminal point of the blue paths.   
This description is technically simpler, and will be used to prove convergence with respect to the $J_1$ metric in \S\ref{sec:tightness}.

\section{Explicit construction and Abelian property}
\label{sec:construction}

In this section we give a more formal definition of the ARW dynamics, and briefly recall the site-wise construction. For details, see~\cite{Rolla19}.

\subsection*{Notation and ARW dynamics}

Let $\N = \{0,1,2,\dots\}$ and $\N_{\sleep} = \N_0 \cup \{\sleep\}$, where $\sleep$ represents a sleeping particle.
For convenience we define $|\sleep|=1$, and $|n|= n$ for $n\in \N$, and write $0<\sleep<1<2<\cdots$.
Also define $\sleep+1=2$ and $n \cdot \sleep = n$ for $n \geq 2$ and $\sleep$ if $n=1$.


The ARW dynamics $ (\eta_t)_t $ is defined as follows.
A site $x$ is \emph{unstable} if it has active particles, i.e., if $\eta_t(x) \geqslant 1$.
At each unstable site $x$, a clock rings at rate $(1+\lambda) \, |\eta_t(x)|$.
When this clock rings, site $x$ is \emph{toppled}, which means that the system goes through the transition $\eta\to \toppling_{x\sleep}\eta$ with probability $\frac{\lambda}{1+\lambda}$, otherwise $\eta\to \toppling_{x\jump}\eta$ with probability $\frac{1}{1+\lambda}$.
These transitions are given by
\[
 \toppling_{x\jump}\eta(z) =
 \begin{cases}
 \eta(x)-1, & z=x, \\
 \eta(y)+1, & z=x+1, \\
 \eta(z), & \mbox{otherwise,}
 \end{cases}
\qquad
 \toppling_{x\sleep}\eta(z) =
 \begin{cases}
 \eta(x) \cdot \sleep, & z=x, \\
 \eta(z), & \mbox{otherwise},
 \end{cases}
\]
and only occur when $\eta(x) \geq 1$.
The operator $\toppling_{x\sleep}$ represents a particle at $x$ trying to fall asleep, which effectively happen if there are no other particles present at $x$.
Otherwise, by definition of $n \cdot \sleep$,
the system state does not change.
The operator $\toppling_{x\jump}$ represents a particle jumping from $x$ to $x+1$, where possible activation of a sleeping particle previously found at $x+1$ is represented by the convention that $\sleep+1=2$.

\subsection*{Site-wise representation and Abelian property}

We now define a field of \emph{instructions} to be read by the active particles. The instructions $\mathcal{I}=({\toppling}^{x,j})_{x\in\Z,j\in\N}$ are i.i.d.\ with $\Pb[ {\toppling}^{x,j}=\toppling_{x\jump} ]=\frac{1}{1+\lambda}$ and $\Pb[ {\toppling}^{x,j}={\toppling_{x\sleep}} ]=\frac{\lambda}{1+\lambda}$.
Using a field of instructions, the operation of \emph{toppling} a site $x$ consists in applying the first instruction available at $x$, and discarding it so that the next unused instruction at $x$ becomes available.
Toppling a site is \emph{legal} if it is unstable.




The ARW dynamics can be recovered from the initial configuration and instructions as follows.
Suppose that every active particle carries a clock which rings according to a Poisson process. Different particles carry independent clocks. When a clock rings for some particle, we topple the site where it is located.
For a system with finite initial configuration, the process obtained this way has the same distribution as the one described above.

The \emph{Abelian property} reads as follows.
For each finite set $V \subseteq \Z$ and initial configuration, if two legal sequences of topplings are contained in $V$ and make each site in $V$ stable, then the resulting configuration is the same.
Using the Abelian property, we can answer many questions of the ARW model by choosing the order in which the sites topple instead of using the Poisson clocks.

\section{Sequential stabilization}
\label{section:seq}

Let $L \in \N$ and, for $x=-L,\dots,0$, let $N_L(x)$ denote the number of particles which jump from $x$ (to $x+1$) when $[-L,0]$ is stabilized.
Moreover, due to the Abelian property, one has that
\[
C_L=N_L(0),
\]
in the sense that $(C_n)_{n \in \N}$ will have the same distribution as the process defined in \S\ref{sec:intro}.
Indeed, we can first stabilize $V_0$, then $V_1$, and so on until $V_L$.
Then writing $k$ for the number of times that $0$ is toppled, $N_L(0)$ will be given by the number of $\toppling_{0\jump}$ found among $\toppling^{0,1},\dots,\toppling^{0,k}$.

In what follows, we will decompose $N_L(x)$ as the sum of two processes, one accounting for the randomness of the initial configuration and another for the randomness of the sleeping instructions. Recall that $\tilde{B}^{+,i}$ is a reflected Brownian motion with diffusion coefficient $r=\sqrt{\sigma_s^2+\sigma_p^2}$.
\begin{proposition}
\label{prop:singletimeconvergence}
For all $i=0,\dots,k$ we have that
\begin{equation}
\label{eq:onlyone}
\left(\varepsilon \, N_{\lfloor \varepsilon^{-2}x_i \rfloor}(\lfloor\varepsilon^{-2}x \rfloor)\right)_{x\geq-x_i} \stackrel{\dd}{\to} (\tilde{B}^{+,i}_x)_{x\geq-x_i}
\end{equation}
as $\varepsilon\to0$ in the metric of uniform convergence on compact intervals of time.
\end{proposition}

For the proof of the proposition above, we will need to describe the sequential stabilization as a reflected random walk.
We will show that for each $L\in\N$ the process $(N_{L}(z))_{z\geq -L}$ is distributed as a random walk started at zero and reflected at zero.

We stabilize $\eta$ on $[-L,0]$ as follows.
Topple site $z=-L$ until it is stable, and denote by $Y_{L}(-L)$ the indicator of the event that the last particle remains passive on $z=-L$.
In case $\eta(-L)=0$, sample $Y_{L}(-L)$ independently of anything else.
By the Abelian property, we have that
\[N_{L}(-L):=[\eta(-L)-Y_{L}(-L)]^+.\]

Note that, after stabilizing $z=-L$, there are $N_{L}(-L)+\eta(-L+1)$ particles at $z=-L+1$.
Now topple site $z=-L+1$ until it is stable, and denote by $Y_{L}(-L+1)$ the indicator of the event that the last particle remains passive on $z=-L+1$.
If there are no particles in $z$, sample $Y_L(z)$ independent of everything else.
Continue this procedure for $z=-L+2,\dots,0$.
Let
\[
T_L(x)=\sum_{y=-L}^x \left(\eta(y)-Y_{L}(y)\right)
\]
and observe that
\[
N_{L}(x)=T_L(x)-\inf_{y=-L,\dots,x} T_L(y).
\]

Write $L=L(\varepsilon)=\lfloor \varepsilon^{-2} x_i \rfloor$.
Note that $Y_L(-L),\dots,Y_L(0),\eta(-L),\dots,\eta(0)$ are independent.
Hence, the increments $(\eta(x)-Y_L(x))_{x\geq -L}$ of $T_L$ are i.i.d.
Since we are assuming that $\zeta=\frac{\lambda}{1+\lambda}$ and $\sigma_p < \infty$, each term in the sum has mean zero and finite variance, so it follows by Donsker's invariance principle that
\begin{equation}
\left( \varepsilon \, T_{L}( \varepsilon^{-2} x ) \right)_{x\in[-x_i,0]}
\stackrel{\dd}{\to}
(\tilde{B}^i_x)_{x\in[-x_i,0]}
\end{equation}
as $\varepsilon\to0$ in the uniform metric.
Since the reflection map $\tilde{B}^i \mapsto \tilde{B}^{+,i}$ in~\eqref{eq:reflect} is continuous, the above convergence implies~\eqref{eq:onlyone}.

\section{Convergence of finite-dimensional projections}
\label{sec:convergence}

Let $k\in\N_0$ and $-x_k\leq -x_{k-1} \leq \dots< -x_0 \leq0$. Recall the definition of $(B^{+,0}_0,{B}^{+,1}_0,\dots,{B}^{+,k}_0)$ in \S\ref{section:limit}.
In this section we prove the following.
\begin{theorem}
[Finite-dimensional convergence]
\label{thm:finitedimensionalconvergence}
We have
\[
\left(\varepsilon \, C_{\lfloor\varepsilon^{-2}x_0\rfloor},\varepsilon \, C_{\lfloor\varepsilon^{-2}x_1\rfloor},\dots,\varepsilon \, C_{\lfloor\varepsilon^{-2}x_k\rfloor}\right)
\stackrel{\dd}{\to}
(B^{+,0}_0,{B}^{+,1}_0,\dots,{B}^{+,k}_0),
\]
as $\epsilon \to 0$.\end{theorem}

Recalling the construction in the previous section, since $C_{\lfloor\varepsilon^{-2}x_i\rfloor} = N_{\lfloor \varepsilon^{-2} x_i \rfloor}(0)$, it is certainly enough to show joint convergence of the counting processes
\[
\left( \left( \vep N_{\lfloor \varepsilon^{-2} x_i \rfloor}(\lfloor\varepsilon^{-2} x\rfloor) \right)_{x \in [-x_i,0]} \right)_{i=0,\dots,k}
\stackrel{\dd}{\to}
\left( \left( B^{+,i}_x \right)_{x \in [-x_i,0]} \right)_{i=0,\dots,k}
.
\]

To keep exposition simpler, we consider the case $k=1$.
There are no differences when considering larger $k$ except for more cluttered notation.
For $x\leq 0$, we define
\[ P(x):=\sum_{y=x}^{0}\left( \zeta-\eta(y)\right)
\]
Given $L\in\N$ and $x\geq -L$, let
\[
S_{L}(x):=\sum_{y=-L}^{x}\left(Y_L(y) - \zeta\right),
\]
where $Y_L(y)$ is as in \S\ref{section:seq}. Actually, we can construct the processes $P,S_L$ and $Y_L$ jointly for different $L$'s. That is, for $-L_1<-L_0<0$,
define 
\begin{equation}\label{eq:2summands}
T_{L_i}(x)=\sum_{y=-L_i}^x \left(\eta(y)-Y_{L_i}(y)\right)
=
P(x)-P(-L_i)-S_{L_i}(x)
\quad i=0,1
\end{equation}
and
\[N_{L_i}(x)=T_{L_i}(x)-\min_{y\in[-L_i,x]} T_{L_i}(y)\quad i=0,1.\]
The key observation is that the $\eta$ terms in~\eqref{eq:2summands} are common for $i=0,1$, whereas the $Y$ terms are independent until $\theta:=\inf\{y:N_{L_0}(y)=N_{L_1}(y)\}$.
After $\theta$, the $Y$ terms are also common for $i=0,1$.

We will need a modified version of the process $T_{L_0}$, whose $S$-component remains independent of $T_{L_1}$ even after $\theta$.
Let
\[
\tilde{Y}_{L_0}(y):=\begin{cases}
Y_{L_0}(y) \quad & \text{ if } y\leq \theta,\\
\tilde{s}(y) \quad & \text{ if } y>\theta,
\end{cases}
\]
where $\tilde{s}$ are i.i.d.\ Bernoulli random variables with parameter $\zeta$, independent of everything else.
Let
\[
\tilde{S}_{L_0}(x):=\sum_{i=-L_0}^x\left(\tilde{Y}_{L_0}(x) - \zeta\right)
\]
and
\[
\tilde{T}_{L_0}(x):=\sum_{i=-L_0}^x \left(\eta(x)-\tilde{Y}_{L_0}(x)\right)
=P(x)- P(-L_0) - \tilde{S}_{L_0}(x).
\]

By Donsker's invariance principle (taking $L_i = \lfloor \varepsilon^{-2} x_i \rfloor$), the triple
\[
\left( (\varepsilon P(\lfloor\varepsilon^{-2}x\rfloor))_{x\leq 0},(\varepsilon S_{\lfloor\varepsilon^{-2} x_1\rfloor}(\lfloor\varepsilon^{-2}x\rfloor))_{x\geq-x_1},(\varepsilon \tilde{S}_{\lfloor\varepsilon^{-2}x_0\rfloor}(\lfloor\varepsilon^{-2}x\rfloor))_{x\geq -x_0} \right)
\]
converges in distribution to three independent Brownian motions
\[
\left( (\mathscr{P}_x)_{x\leq 0},(\mathscr{S}_t^1)_{x\geq -x_1},(\tilde{\mathscr{S}}^0_t)_{x\geq -x_0} \right)
\]
with diffusion coefficients $\sigma_p,\sigma_s,\sigma_s$ respectively and started at 0 (i.e., $\mathscr{P}_{0}=\mathscr{S}^1_{-x_1}=\tilde{\mathscr{S}}^0_{-x_0}=0$).

Let 
\[\mathscr{T}^1_x:=\mathscr{P}_x-\mathscr{P}_{-x_1} - \mathscr{S}^1_x, \quad 0\geq x\geq-x_1\] and
\[\tilde{\mathscr{T}}^0_x:=\mathscr{P}_x-\mathscr{P}_{-x_0} - \tilde{\mathscr{S}}^0_x, \quad x\geq -x_0.\]
Let also 
\[\mathscr{N}^1_x=\mathscr{T}^1_x-\inf_{s\in[-x_1,x]}\mathscr{T}^1_s\] and 
\[\tilde{\mathscr{N}}^0_x=\tilde{\mathscr{T}}^0_x-\inf_{s\in[-x_0,x]}\tilde{\mathscr{T}}^0_s.\]

Define \[N^{\varepsilon,0}_x:=\varepsilon\, N_{\lfloor \varepsilon^{-2} x_0 \rfloor}(\varepsilon^{-2} x),\] also define $\tilde{N}^{\varepsilon,0}_x$, $N^{\varepsilon,1}_x$, $\tilde{T}^{\varepsilon,0}_x$ and ${T}^{\varepsilon,1}_x$ analogously.
As in the proof of Proposition~\ref{prop:singletimeconvergence}, by invariance principle and continuity of the map $(f(t))_{t\geq 0}\mapsto (\inf_{s \leq t}f(s))_{t\geq 0}$ under the uniform metric, we have that
\begin{equation}
\label{eq:a.s.convergence}
\begin{aligned}
((\tilde{T}^{\varepsilon,0}_x)_{x\geq -x_0},
({T}^{\varepsilon,1}_x)_{x\geq -x_1},
(\tilde{N}^{\varepsilon,0}_x)_{x\geq -x_0},(N^{\varepsilon,1}_x)_{x\geq -x_1})\\
\text{converges in distribution, as }\vep\to0\text{ to}\\
((\tilde{\mathscr{T}}^0_x)_{x\geq -x_0},({\mathscr{T}}^1_x)_{x\geq -x_1},(\tilde{\mathscr{N}}^0_x)_{x\geq -x_0},({\mathscr{N}}^1_x)_{x\geq -x_1})
\end{aligned}
\end{equation}
uniformly over compacts. By the Skorohod representation theorem we can (and will) assume that the convergence above holds almost surely.

We still have to show that ${N}^{\varepsilon,0}$ converges to $\mathscr{N}^0$ defined now.
Let
\[
\tau:=\inf\{x\geq -x_0: \tilde{\mathscr{N}}^0_x = \mathscr{N}^1_x \}
\
\text{ and }
\
\mathscr{N}^0_x :=
\begin{cases}
\tilde{\mathscr{N}}^0_x, & x \in[-x_0,\tau),\\
\mathscr{N}^1_x, & x \geq \tau.
\end{cases}
\]

To prove ${N}^{\varepsilon,1} \to \mathscr{N}^1$ we consider the coalescing time of discrete processes.
Writing
\[
\tau^\varepsilon:=\inf\{x\geq -x_0:\tilde{N}^{\varepsilon,0}_x=N^{\varepsilon,1}_x\}
,
\]
we have
\[
{N}^{\varepsilon,0} =
\begin{cases}
\tilde{{N}}^{\varepsilon,0}_x, & x \in[-x_0,\tau),\\
{N}^{\varepsilon,1}_x, & x \geq \tau^\varepsilon,
\end{cases}
\]
so to conclude the proof it suffices to show that $\tau^\varepsilon \to \tau$ a.s.

Since the first time that two paths meet is a lower semi-continuous function of the paths, we have
$\liminf\limits_{\varepsilon\to0} \tau^{\varepsilon}\geq \tau$ a.s. It remains to show that
\begin{equation}\label{eq:cotacara}
\limsup_{\varepsilon\to0} \tau^{\varepsilon} \leq \tau \quad\text{a.s.}
\end{equation}

Since $N_{\lfloor \epsilon^{-2}x_1\rfloor}(x) - \tilde{N}_{\lfloor \epsilon^{-2}x_0\rfloor}(x) \geq 0$ at $x=-x_0$ and this difference only jumps by $0$ or $\pm 1$, it suffices to show the following claim:
\begin{equation}
\label{eq:claim}
\text{Given any $\delta>0$, a.s., for all $\vep$ small enough, there is $x \le \tau+\delta$ with
$N^{\varepsilon,1}_x \leq \tilde{N}^{\varepsilon,0}_x$.}
\end{equation}

Since $\tau$ is a stopping time for the pair $(\mathscr{T}^1,\mathscr{T}^0)$, by the strong Markov property the processes $(\Delta_x)_{x\geq 0} := (\mathscr{T}^1_{\tau+x}-\mathscr{T}^1_\tau)_{x\geq0}$ and $(\tilde\Delta_x)_{x\geq 0} := (\tilde{\mathscr{T}}^0_{\tau+x} - \mathscr{T}^1_{\tau+x})_{x\geq0}$ are distributed as Brownian motions started at value $0$.
Hence, a.s.~there is a point $z^* \in [0,\delta]$ such that $\tilde{\Delta}_{z^*} > 0$.
Moreover, since $\dd \mathscr{N} \geq \dd \mathscr{T}$ due to reflection, we either have ${\mathscr{N}}^1_{\tau+z^*}=
{\mathscr{N}}^1_\tau + \Delta_{z^*}$ or ${\mathscr{N}}^1_{\tau+z^*}>{\mathscr{N}}^1_\tau + \Delta_{z^*}$. We will distinguish between those two cases.
 In the first case
\begin{equation}
\label{eq:display}
{\mathscr{N}}^1_{\tau+z^*}
=
{\mathscr{N}}^1_\tau + \Delta_{z^*}
=
\tilde{\mathscr{N}}^0_\tau + \Delta_{z^*}
<
\tilde{\mathscr{N}}^0_{\tau}
+ \Delta_{z^*}
+ \tilde{\Delta}_{z^*}
\leq
\tilde{\mathscr{N}}^0_{\tau+z^*}.
\end{equation}
In particular, we have obtained the strict inequality 
\[{\mathscr{N}}^1_{\tau+z^\ast} < \tilde{\mathscr{N}}^0_{\tau+z^\ast}.\]
By this inequality and~\eqref{eq:a.s.convergence}, it follows that, for $n$ large enough, 
\[N^{\vep,1}_{\tau+z^\ast} < \tilde{N}^{\vep,0}_{\tau+z^\ast}.\]
In this first case,~\eqref{eq:claim} follows directly from this inequality.

The case 
\[{\mathscr{N}}^1_{\tau+z^\ast}>{\mathscr{N}}^1_\tau + \Delta_{z^\ast}\] is subtler. We work it by observing that the above inequality is equivalent to
\[{\mathscr{T}}^1_{\tau+z^\ast}-\min_{x\in[-x_1,\tau+z^\ast]}{\mathscr{T}}^1_{x}>{\mathscr{T}}^1_{\tau}-\min_{x\in[-x_1,\tau]}{\mathscr{T}}^1_{x}+ \mathscr{T}^1_{\tau+z^\ast}-\mathscr{T}^1_\tau,\]
which, in turn, is equivalent to
\[
\min_{x\in[-x_1,\tau+z^*]} \mathscr{T}^1_x< \min_{x\in[-x_1, \tau^*]} \mathscr{T}^1_x.
\]
From~\eqref{eq:a.s.convergence}, this implies that, for $\vep$ small enough,
\[
\min_{x\in[-x_1,\tau+z^*]}{T}^{\varepsilon,1}_{x}
<
\min_{x\in[-x_1,\tau]}{T}^{\varepsilon,1}_{x}
,
\]
in which case there is $x^*\in[\tau,\tau+z^*]$ such that
\[
{T}^{\varepsilon,1}_{x^*} =
\min_{s\in[-x_1,x^*]}{T}^{\varepsilon,1}_{s}
,
\]
whence
\[
{N}^{\varepsilon,1}_{x^*} = 0 \leq 
\tilde{N}^{\varepsilon,0}_{x^*}
.
\]
Hence~\eqref{eq:claim} holds, and this completes the proof of Theorem~\ref{thm:finitedimensionalconvergence}.

\section{Convergence in the $J_1$ metric}
\label{sec:tightness}

\begin{figure}[b]
\centering
\includegraphics[width=\textwidth]{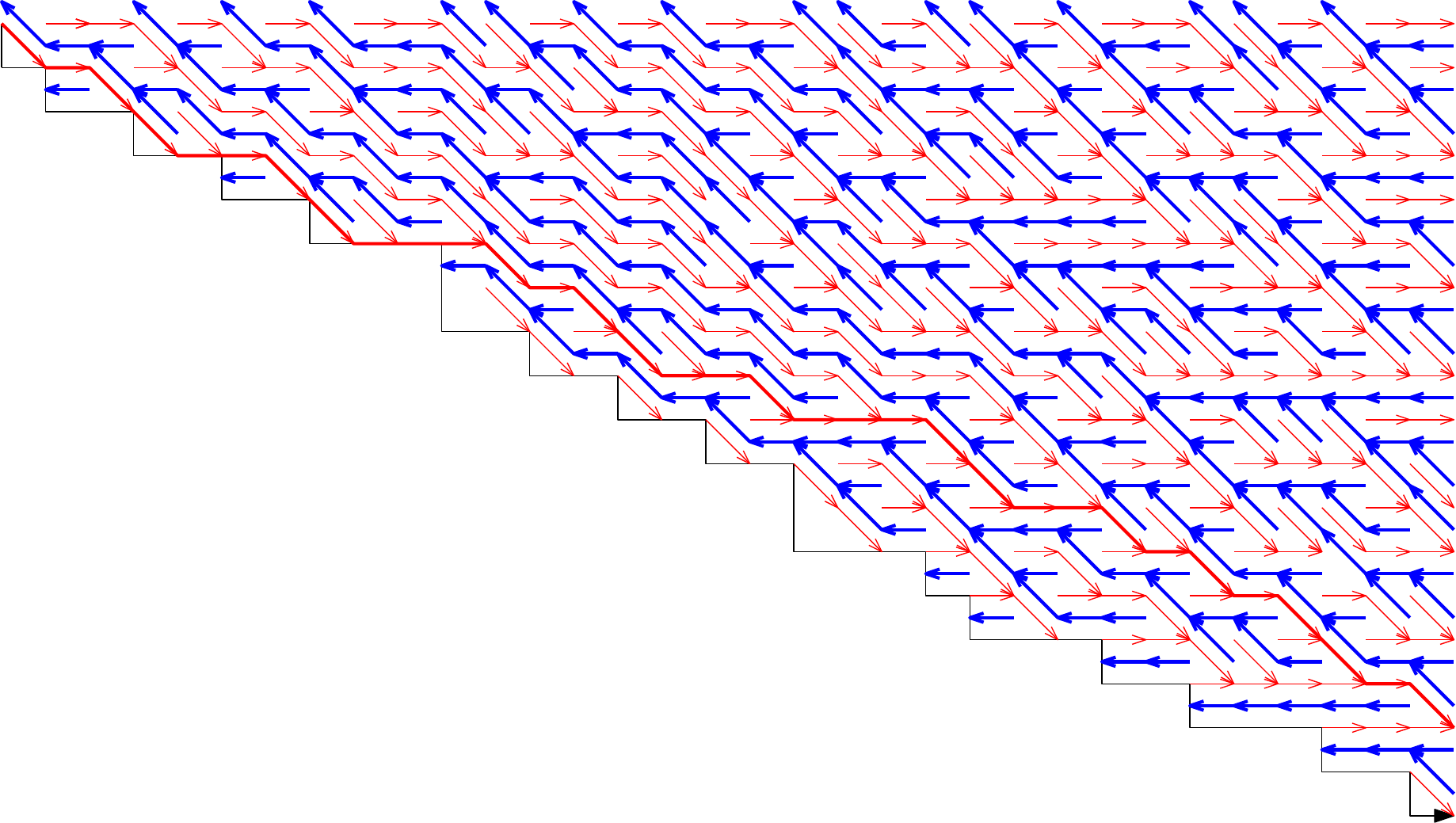}
\caption{Black path, red paths and dual blue paths.}
\label{fig:dualpaths}
\end{figure}
 In this section we prove the following.

\begin{theorem}
\label{thm:tight}
The family of processes $(\vep C_{\lfloor \vep^{-2}x\rfloor})_{x \geq 0}$ indexed by $\varepsilon \in (0,1]$ is tight in the Skorohod $J_1$ metric.
\end{theorem}

\begin{figure}[b]
\centering
\includegraphics[width=\textwidth]{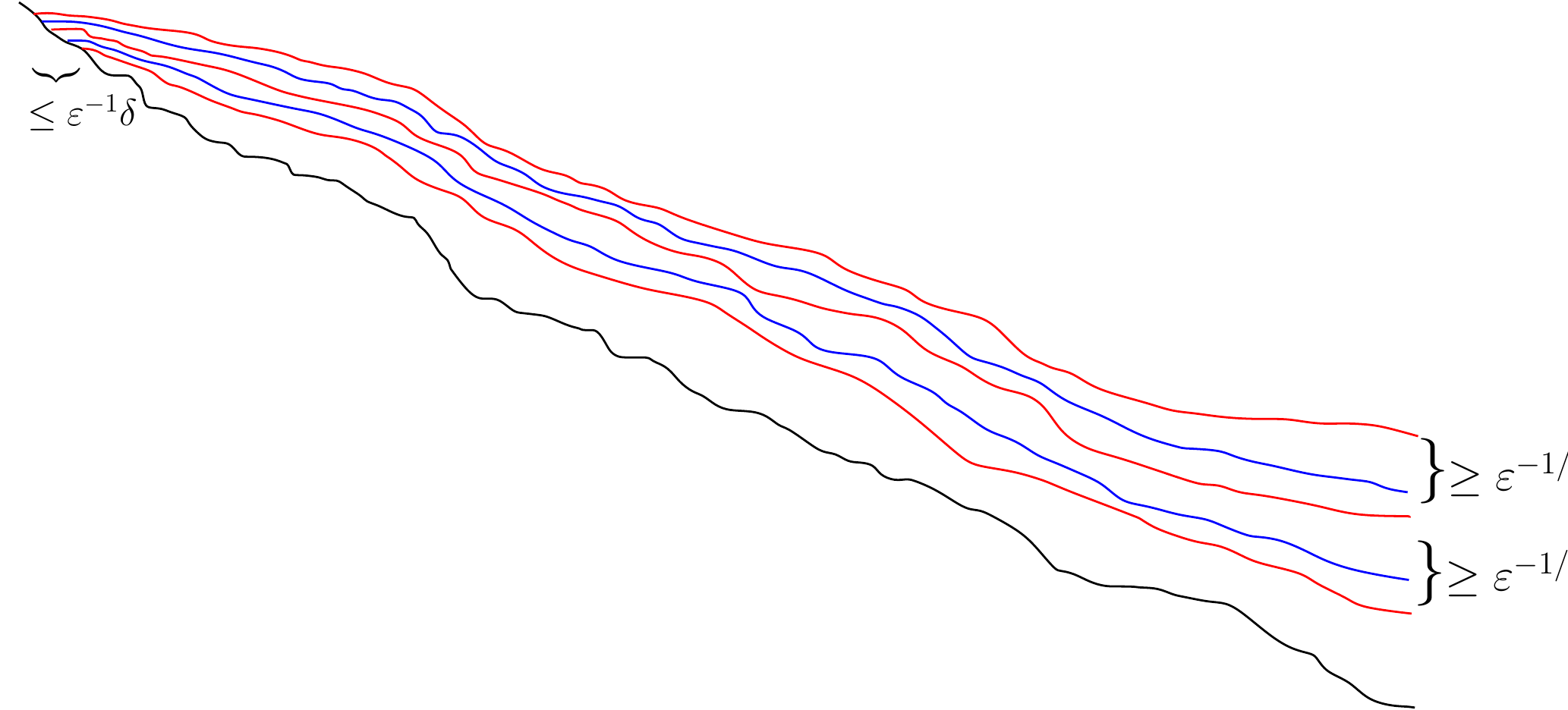}
\caption{Event $\mathfrak{B}(\varepsilon,\delta,\gamma,M,K)$ in terms of black, red and blue paths.}
\label{fig:event}
\end{figure}

Since the processes $(\varepsilon\, C_{\lfloor \varepsilon^{-2} x \rfloor})_{x\geq 0}$ are non-decreasing for each $\vep$, by Theorem~\ref{thm:finitedimensionalconvergence} they converge in the $M_1$-metric to a process $(\mathscr{C}_x)_{x \geq 0}$.
Moreover, in order to get tightness in the $J_1$ metric, it suffices (see~\cite[Theorem 13.2]{Billingsley99}) to show that, for each $K < \infty$ and $\gamma>0$ fixed,
\begin{equation}
\label{eq:boundedtight}
\lim_{R \to \infty} \limsup_{\varepsilon \to 0} \Pb\left[\sup_{x\in[0,K]} | C_x^\varepsilon |\geq R\right]=0
\end{equation}
and
\begin{equation}
\label{eq:twojumps}
\lim_{\delta \to 0} \limsup_{\varepsilon\to0}\Pb\left[\sup| C^{\varepsilon}_{x_1}-C^{\varepsilon}_{x_0}| \wedge | C^\varepsilon_{x_1}-C^\varepsilon_{x_2}|\geq \gamma\right]=0
.
\end{equation}
where the last supremum is taken over triples $x_0\leq x_1\leq x_2$ which satisfy
$x_2\leq K$ and $x_2-x_0\leq \delta$.
Note that~\eqref{eq:boundedtight} follows directly from Theorem~\ref{thm:finitedimensionalconvergence}.
We only have to establish~\eqref{eq:twojumps}.
The event whose probability we want to control is that of having two macroscopic jumps in a short interval.

For this proof we will use a graphical construction for the flow process, which we now proceed to describe.
First, consider the \emph{black path} $B:(-\infty,0]\to\Z$, given by $B(x)=\sum_{y=x}^0\eta(y)$, where we recall that $\eta(y)$ is the initial number of particles at $y$. As in Figure~\ref{fig:dualpaths}, for every lattice point $(x,y) \in (-\infty,0]\times\N$ above the graph of $B$, place a red arrow which either points to $(x+1,y-1)$ or to $(x+1,y)$, independently with probability $\frac{\lambda}{1+\lambda}=\zeta$ and $\frac{1}{1+\lambda}=1-\zeta$, respectively. Let $x\in(-\infty,0]$. We denote by $R_x:[x,0]\to \Z$ the path which starts at $(x,B(x))$, follows the red arrows and is reflected on the black path $B$. The resulting path will be referred as the \emph{red path}. In Figure~\ref{fig:dualpaths}, the red path is depicted in bold red. Note that, by associating the red arrows which point down with sleeping instructions and the horizontal red arrows with no-sleep instructions, we see that $R_x(0)$ is distributed as $N_{-x}(0)=C_{x}$.
Moreover, one can use the same black path and red arrows to get any joint distribution $(C_{x_0},\dots,C_{x_k})$, $-x_k \leq -x_{k-1} \leq \cdots \leq -x_0 \leq 0$: simply consider the collection of red paths $R_{-x_k},\dots,R_{-x_0}$ constructed using the same black path and the same red arrows. The random vector $(R_{-x_0}(0),\dots,R_{-x_k}(0))$ will have the distribution of $(C_{x_0},\dots,C_{x_k})$. Observe that the red paths are independent until they meet, after which they coalesce.

Now, consider a collection of \emph{blue arrows} which is dual to the red arrows. That is, for every $(x,y)\in \Z^2$ above the graph of the black path, there is a (backwards) blue arrow emanating from $(x,y-\frac{1}{2})$.
That blue arrow points to $(x-1,y+\frac{1}{2})$ if the red arrow staring at $(x-1,y)$ points to $(x,y-1)$.
Otherwise, the blue arrows points to $(x-1,y-\frac{1}{2})$ (See Figure~\ref{fig:dualpaths}). By following the blue arrows we can construct a collection of \emph{blue paths} which are dual to the red ones. The blue paths should be read from right to left.
The blue paths jump at each time with probability $\zeta$. Also, the blue paths are killed when they meet the black path.
For any $y\in\Z$, let $A_y$ be the blue path starting at $(0,y+\frac{1}{2})$. Observe that, by construction, the blue paths cannot intersect the red paths.

 In the following reasoning, we will consider jointly two different blue paths $A_{y_1}, A_{y_2}$. For convenience, we will consider a modified version of the processes $\tilde{A}_{y_1},\tilde{A}_{y_2}$ such that they evolve independently after collision (do not coalesce). Also, they don't get killed at intersecting the black path, they evolve independently. We do this because now $\tilde{A}_{y_1},\tilde{A}_{y_2}$ and $B$ are three independent random walks and we can apply Donsker's invariance principle.
 
 We consider the (modified) blue path $\tilde{A}_y$ as a stochastic process $(\tilde{A}_y(x))_{x\leq0}$, where, for any $x\leq 0$, $\tilde{A}_y(x)$ is the position at $x$ of the (modified) blue path started at $(0,y+\frac{1}{2})$.

Let \[(\tilde{A}^\varepsilon_{y_1}(x))_{x\leq 0}:=(\varepsilon \tilde{A}_{\lfloor \vep^{-1} y_1 \rfloor}(\lfloor\varepsilon^{-2}x\rfloor) - \zeta x)_{x\leq0}\]
\[(\tilde{A}^\varepsilon_{y_2}(x))_{x\leq 0}:=(\varepsilon \tilde{A}_{\lfloor \vep^{-1} y_2 \rfloor}(\lfloor \varepsilon^{-2}x \rfloor) - \zeta x)_{x\leq 0}\]
\[(B^\varepsilon(x))_{x\leq0}:=(\varepsilon B(\lfloor\varepsilon^{-2}x\rfloor) - \zeta x)_{x\leq0}.\]
By Donsker's invariance principle, we get that
\begin{equation}\label{eq:threeind}
((\tilde{A}^\vep_{y_1}(x))_{x\leq0},(\tilde{A}^\vep_{y_2}(x))_{x\leq0},(\tilde{B}^\vep(x))_{x\leq0})\to
((\tilde{\mathcal{A}}_{y_1}(x))_{x\leq0},(\tilde{\mathcal{A}}_{y_2}(x))_{x\leq0},(\mathcal{B}(x))_{x\leq0}),
\end{equation}
where $\tilde{\mathcal{A}}_{y_1},\tilde{\mathcal{A}}_{y_2},\mathcal{B}$ are three (time-reversed) independent Brownian motions started at $y_1,y_2$ and $0$ respectively.
By the Skorohod Representation theorem, we can (and will) assume that the convergence holds almost surely.

It follows that, under the event in~\eqref{eq:twojumps}, both $\tilde{A}^\varepsilon_{y_1}$ and $\tilde{A}^\varepsilon_{y_2}$ (when read from $0$ to $-\infty$) intersect $B^\varepsilon$ for the first time in a time window of length smaller than $\delta$ (see Figure~\ref{fig:event}). 
By~\eqref{eq:threeind}, the probability of the event above converges to the probability that two independent Brownian motions $\tilde{\mathcal{A}}_{y_1},\tilde{\mathcal{A}}_{y_2}$ intersect a third B.m., $\mathcal{B}$, also independent, in a time window of length smaller than $\delta$.
Using the continuity of the Brownian motion, we get that, as $\delta\to 0$, this converges to the probability that three independent Brownian motions eventually meet at the same point at the same time.
The latter probability is $0$ and this will finish the proof.
To complete the proof, we have to explain how to choose the initial points of the blue paths $y_1,y_2$.

We will choose $M$ large but fixed.
It suffices to deal with the case $C_{K}^{\varepsilon} \geq M$, since Theorem~\ref{thm:finitedimensionalconvergence} readily implies that
\[\lim_{M\to\infty}\limsup_{\epsilon\to0} \Pb\left[C_{K}^\varepsilon> M\right]=0.\]

Therefore, it is enough to show that, for every $M>0$ and $\gamma>0$,
\begin{equation}
\lim_{\delta \to 0} \limsup_{\varepsilon\to0}\Pb\left[\sup_{\stackrel{x_0\leq x_1\leq x_2\in[0,K]}{x_2-x_0\leq \delta}} (C^{\varepsilon}_{x_1}-C^{\varepsilon}_{ x_0}) \wedge (C^{\varepsilon}_{x_2}-C^{\varepsilon}_{x_1})\geq \gamma;C_{x_2}^\varepsilon\leq M\right]=0.
\end{equation}
Let $\mathfrak{B}(\varepsilon,\delta,\gamma,M,K)$ be the event inside the above probability (see Figure~\ref{fig:event}).
If $\mathfrak{B}(\varepsilon,\delta,\gamma,M,K)$ holds, then, there exists $x_0^*<x_1^*<x_2^*$ with $x_2^*-x_0^*\leq \delta$ such that both $C^{\varepsilon}_{x_1^*}-C^{\varepsilon}_{x_0^*}$ and $C^{\varepsilon}_{x_2^*}-C^{\varepsilon}_{x_1^*}$ are greater than $\gamma$. Moreover, $C^\vep_{x_2^*}\leq M$. 
 
 Now, we partition $[0,M]$ into intervals of size $\frac{\gamma}{2}$.
 We consider the intervals $I_i=[\frac{\gamma i}{2},\frac{\gamma (i+1)}{2}]$, $i=0,\dots,\lceil \frac{2M}{\gamma}\rceil$.
 By the discussion on the paragraph above, there must be indices $i_1<i_2\in\{0,\dots,\lceil \frac{2M}{\gamma}\rceil\}$ such that, $I_{i_1}$ is contained in $[C^{\varepsilon}_{x_0^*},C^{\varepsilon}_{x_1^*}]$ and $I_{i_2}$ is contained in $[C^{\varepsilon}_{x_1^*},C^{\varepsilon}_{x_2^*}]$. Let $\mathfrak{B}^{\varepsilon,\delta}_{i_1,i_2}$ be the event described just above.

 Then, 
 \begin{equation}\label{eq:split}
\begin{aligned}
\mathbb{P}[\mathfrak{B}(\varepsilon,\delta,\gamma,M,K)]&\leq\sum_{i_1< i_2\in\{0,\dots,\lceil\frac{2M}{\gamma}\rceil\}}
\mathbb{P}\left [ {\mathfrak B}^{\varepsilon,\delta}_{i_1,i_2}\right]\\
&\leq \left(\left\lceil\frac{2M}{\gamma}\right\rceil+1\right)^2\sup_{i_1< i_2\in\left\{0,\dots,\left\lceil\frac{2M}{\gamma}\right\rceil\right\}}\mathbb{P}\left [ {\mathfrak B}^{\varepsilon,\delta}_{i_1,i_2}\right].
\end{aligned}
\end{equation}
 
We now bound the probability of ${\mathfrak B}^{\varepsilon,\delta}_{i_1,i_2}$ for $i_1< i_2$. Since blue paths do not intersect red paths under $\mathfrak{B}^{\varepsilon,\delta}_{i,j}$, any blue path started between $C^\vep_{x_0}$ and $C^\vep_{x_1}$ will remain between $C^\vep_{x_0}$ and $C^\vep_{x_1}$ and therefore will intersect $B^\vep$ at some point $z_1^*\in[-x_1,-x_0]$. Hence, since $I_{i_1}$ is contained in $[C^{\varepsilon}_{x_0},C^{\varepsilon}_{x_1}]$, we have that $\tilde{A}^\vep_{\frac{\gamma}{2} i_1}$ intersects $B^\vep$ in $[-x_1,-x_0]$. Analogously, $\tilde{A}^\vep_{\frac{\gamma }{2}i_2}$ intersects $B^\vep$ in $[-x_2,-x_1]$. Hence, since $x_2-x_1\leq \delta$, both $\tilde{A}^\vep_{\frac{\gamma}{2} i_1}$ and $\tilde{A}^\vep_{\frac{\gamma}{2} i_2}$ intersect $B^\varepsilon$ in a time interval of size $\delta$.

That is, let $\tau^\varepsilon_i:=\inf\{t\geq0: \tilde{A}^\varepsilon_{\frac{\gamma}{2}i}(t)=B^\varepsilon(t)\}.$ 
Then 
\[
\mathfrak{B}^{\varepsilon,\delta}_{i_1,i_2}\subseteq\{|\tau^\varepsilon_{i_1}-\tau^\varepsilon_{i_2}|\leq \delta\}.
\]
 
Therefore, by~\eqref{eq:threeind},
 \[\limsup_{\varepsilon\to0}\mathbb{P}\left [ {\mathfrak B}^{\varepsilon,\delta}_{i_1,i_2}\right]\leq\mathbb{P}\left[|\tau_{i_1}-\tau_{i_2}|\leq \delta\right]\] 
 where $\tau_i:=\inf_{t\geq 0}\{ \mathcal{A}_{\frac{\gamma}{2}i}(t)=\mathcal{B}(t)\}$.
 
 Hence,
 \[
 \lim_{\delta\to0}\limsup_{\vep\to0}\mathbb{P}[{\mathfrak B}^{\varepsilon,\delta}_{i_1,i_2}]\leq \lim_{\delta\to0}\mathbb{P}[|\tau_{i_1}-\tau_{i_2}|\leq \delta]=\mathbb{P}[\tau_{i_1}=\tau_{i_2}].
 \]

The proof is finished by noticing that $\tau_{i_1}=\tau_{i_2}$ implies that the three-dimensional Brownian motion $(\mathcal{A}_{\frac{\gamma}{2}i_1},\mathcal{A}_{\frac{\gamma}{2}i_2},\mathcal{B})$ intersects the line $\{(x,y,z)\in\mathbb{R}^3:x=y=z\}$, and that event has zero probability.
 
 
\section{The scaling limit is a pure-jump process}
\label{sec:purejump}

In this section we prove Theorem~\ref{thm:purejump}.
By Theorem~\ref{thm:convergence}, the non-decreasing càdlag process $\left( \CC^\rho_{x}\right)_{x \ge 0}$ is well-defined.
Following the description of the previous section, it can be constructed directly from a black Brownian motion $\mathcal{B}=(\mathcal{B}_x)_{x \le 0}$ started at $\mathcal{B}_0 = 0$ and blue coalescing Brownian motions $(\mathcal{A}_y(x))_{x \le 0}$ starting from $\mathcal{A}_y(0)=y$ and independent of $\mathcal{B}$, having diffusion coefficients $1$ and $\rho$ respectively.


The construction is as follows.
First observe that the red paths can be recovered from the black and blue paths, so they are not needed in the construction.
Moreover, blue paths are independent of each other until they coalesce, and independent of the black path until they are killed by it.
Furthermore, since blue paths coalesce, killing them upon meeting the black path is irrelevant and we can disconsider it.
In the scaling limit, this collection of blue paths converges to the a family of paths~\cite{Arratia79,FontesIsopiNewmanRavishankar04} consisting of independent coalescing paths $(\mathcal{A}_y(x))_{x \le 0}$ indexed by $y>0$, each one started from $\mathcal{A}_y(0)=y$.

This family satisfies the following.
Let $\T_y = \inf\{x \ge 0: \mathcal{A}_y(-x)=\mathcal{B}_{-x}\}$.
Then a.s.\ $0 < \T_y < \infty$ for every $y>0$, $y\mapsto \T_y$, is non-decreasing, the set of values $\{\T_y : y\in[\delta,K]\}$ is finite for every $0<\delta<K<\infty$, $\lim\limits_{y\to 0^+} \T_y =0$, and $\lim\limits_{y \to \infty} \T_y = \infty$.

To be self-contained, let us justify the statements in the previous paragraph more carefully.
Write $\mathbb{Q}_+^*=\{y_n\}_{n\in\N}$.
Take $\mathcal{A}_{y_1}=(\mathcal{A}_{y_1}(x))_{x\le 0}$ starting from $\mathcal{A}_{y_1}(0)=y_1$.
For each $n$, take $\mathcal{A}_{y_n}=(\mathcal{A}_{y_n}(x))_{x\le 0}$ starting from $\mathcal{A}_{y_n}(0)=y_n$, independent of $\mathcal{A}_{y_1},\dots,\mathcal{A}_{y_{n-1}}$ until the first point (that is, highest $x$) where it meets one of them, and equal to that one after such time (that is, for lower values of $x$).
Now for $y>0$ rational, let $\T_y = \inf\{x \ge 0: \mathcal{A}_y(-x)=\mathcal{B}_{-x}\}$.
Then a.s.\ $0<\T_y<\infty$ for every $y$.
By coalescence, $\mathcal{A}_y \le \mathcal{A}_{y'}$ for $y<y'$, hence $y\mapsto \T_y$, is non-decreasing.
Furthermore, using Borel-Cantelli one can show that
$\lim_{y\to 0^+} \T_y =0$ and $\lim_{y \to \infty} \T_y = \infty$.
Finally, by well-known properties of coalescing Brownian motions~\cite{Arratia79,FontesIsopiNewmanRavishankar04}, for each $0<a<b<\infty$ and $\vep>0$, the set $\{\mathcal{A}_y(-\vep):a<y<b\}$ is a.s.\ finite.
As a consequence of the two last properties, the set $\{\T_y : y\in[\delta,K]\}$ is finite for every $0<\delta<K<\infty$.

To conclude, following the description of the previous section, we can define
\[
\CC^\rho_x := 
\inf \left\{ y > 0 : \T_y > x \right\}
,\quad
x \ge 0.
\]
By the remarks of the previous paragraph, a.s.\ for every $0<a<b<\infty$ the process $\left( \CC^\rho_{x}\right)_{x \in [a,b]}$ takes only finitely many values, and therefore it is a pure-jump process.

\section*{Acknowledgment}

We thank V.~Sidoravicius and J.~Yu for helpful discussions.
Manuel Cabezas was supported by Iniciativa Cient\'{i}fica Milenio NC120062.

\parskip 0pt
\setstretch{1}
\small\bibliographystyle{bib/leoabbrv}
\bibliography{bib/leo}

\begin{thebibliography}{10}
\expandafter\ifx\csname urlstyle\endcsname\relax
  \providecommand{\doi}[1]{doi:\discretionary{}{}{}#1}\else
  \providecommand{\doi}{doi:\discretionary{}{}{}\begingroup
  \urlstyle{rm}\Url}\fi

\bibitem{AmirGurel-Gurevich10}
G.~\textsc{Amir}, O.~\textsc{Gurel-Gurevich}.
\newblock \emph{On fixation of activated random walks}.
\newblock Electron Commun Probab \textbf{15}:119--123, 2010.
\newblock \doi{10.1214/ECP.v15-1536}.

\bibitem{Arratia79}
R.~\textsc{Arratia}.
\newblock \emph{Coalescing {Brownian} Motions on the Line}.
\newblock Ph.D. thesis, University of Wisconsin, Madison, 1979.

\bibitem{AsselahSchapiraRolla19}
A.~\textsc{Asselah}, B.~\textsc{Schapira}, L.~T. \textsc{Rolla}.
\newblock \emph{Diffusive bounds for the critical density of activated random
  walks}, 2019.
\newblock Preprint. \href{http://arxiv.org/abs/1907.12694}{arXiv:1907.12694}.

\bibitem{BasuGangulyHoffman18}
R.~\textsc{Basu}, S.~\textsc{Ganguly}, C.~\textsc{Hoffman}.
\newblock \emph{Non-fixation for conservative stochastic dynamics on the line}.
\newblock Comm Math Phys \textbf{358}:1151--1185, 2018.
\newblock \doi{10.1007/s00220-017-3059-7}.

\bibitem{BasuGangulyHoffmanRichey19}
R.~\textsc{Basu}, S.~\textsc{Ganguly}, C.~\textsc{Hoffman}, J.~\textsc{Richey}.
\newblock \emph{Activated random walk on a cycle}.
\newblock Ann Inst Henri Poincar\'{e} Probab Stat \textbf{55}:1258--1277, 2019.
\newblock \doi{10.1214/18-aihp918}.

\bibitem{Billingsley99}
P.~\textsc{Billingsley}.
\newblock \emph{Convergence of probability measures}.
\newblock Wiley Series in Probability and Statistics: Probability and
  Statistics. John Wiley \& Sons, Inc., New York, 2 edn., 1999.
\newblock \doi{10.1002/9780470316962}.

\bibitem{CabezasRollaSidoravicius14}
M.~\textsc{Cabezas}, L.~T. \textsc{Rolla}, V.~\textsc{Sidoravicius}.
\newblock \emph{Non-equilibrium phase transitions: Activated random walks at
  criticality}.
\newblock J Stat Phys \textbf{155}:1112--1125, 2014.
\newblock \doi{10.1007/s10955-013-0909-3}.

\bibitem{CabezasRollaSidoravicius18}
---{}---{}---.
\newblock \emph{Recurrence and density decay for diffusion-limited annihilating
  systems}.
\newblock Probab Theory Relat Fields \textbf{170}:587--615, 2018.
\newblock \doi{10.1007/s00440-017-0763-3}.

\bibitem{DickmanRollaSidoravicius10}
R.~\textsc{Dickman}, L.~T. \textsc{Rolla}, V.~\textsc{Sidoravicius}.
\newblock \emph{Activated random walkers: Facts, conjectures and challenges}.
\newblock J Stat Phys \textbf{138}:126--142, 2010.
\newblock \doi{10.1007/s10955-009-9918-7}.

\bibitem{FontesIsopiNewmanRavishankar04}
L.~R.~G. \textsc{Fontes}, M.~\textsc{Isopi}, C.~M. \textsc{Newman},
  K.~\textsc{Ravishankar}.
\newblock \emph{The {B}rownian web: characterization and convergence}.
\newblock Ann Probab \textbf{32}:2857--2883, 2004.
\newblock \doi{10.1214/009117904000000568}.

\bibitem{Rolla19}
L.~T. \textsc{Rolla}.
\newblock \emph{Activated random walks on {$Z^d$}}, 2019.
\newblock Preprint. \href{http://arxiv.org/abs/1906.05037}{arXiv:1906.05037}.

\bibitem{RollaSidoravicius12}
L.~T. \textsc{Rolla}, V.~\textsc{Sidoravicius}.
\newblock \emph{Absorbing-state phase transition for driven-dissipative
  stochastic dynamics on {$Z$}}.
\newblock Invent Math \textbf{188}:127--150, 2012.
\newblock \doi{10.1007/s00222-011-0344-5}.

\bibitem{RollaTournier18}
L.~T. \textsc{Rolla}, L.~\textsc{Tournier}.
\newblock \emph{Non-fixation for biased activated random walks}.
\newblock Ann Inst H Poincar{\'e} Probab Statist \textbf{54}:938--951, 2018.
\newblock \doi{10.1214/17-AIHP827}.

\bibitem{Shellef10}
E.~\textsc{Shellef}.
\newblock \emph{Nonfixation for activated random walks}.
\newblock ALEA Lat Am J Probab Math Stat \textbf{7}:137--149, 2010.
\newblock \href{http://alea.impa.br/articles/v7/07-07.pdf}{pdf}.

\bibitem{SidoraviciusTeixeira17}
V.~\textsc{Sidoravicius}, A.~\textsc{Teixeira}.
\newblock \emph{Absorbing-state transition for stochastic sandpiles and
  activated random walks}.
\newblock Electron J Probab \textbf{22}:33, 2017.
\newblock \doi{10.1214/17-EJP50}.

\bibitem{StaufferTaggi18}
A.~\textsc{Stauffer}, L.~\textsc{Taggi}.
\newblock \emph{Critical density of activated random walks on transitive
  graphs}.
\newblock Ann Probab \textbf{46}:2190--2220, 2018.
\newblock \doi{10.1214/17-AOP1224}.

\bibitem{Taggi16}
L.~\textsc{Taggi}.
\newblock \emph{Absorbing-state phase transition in biased activated random
  walk}.
\newblock Electron J Probab \textbf{21}:13, 2016.
\newblock \doi{10.1214/16-EJP4275}.

\bibitem{Taggi19}
---{}---{}---.
\newblock \emph{Active phase for activated random walks on {$\Bbb{Z}^d$},
  {$d\geq3$}, with density less than one and arbitrary sleeping rate}.
\newblock Ann Inst Henri Poincar\'{e} Probab Stat \textbf{55}:1751--1764, 2019.
\newblock \doi{10.1214/18-aihp933}.

\end{thebibliography}

\end{document}